\documentclass[11pt]{article}
\usepackage{mathrsfs}
\usepackage{amscd}
\usepackage{amsmath,amsfonts,amssymb,amscd}
\usepackage{indentfirst,graphics,epsfig,psfrag}
\input{epsf}
\usepackage{ifpdf}
\usepackage{enumerate}
\usepackage{appendix}
\usepackage{enumerate}
\usepackage{lineno}
\usepackage{color}
\makeatletter \@addtoreset{figure}{section} \makeatother
\makeatletter
\long\def\@makecaption#1#2{%
   \vskip 10\p@
   \setbox\@tempboxa\hbox{{#1}\ \ #2}%
   \ifdim \wd\@tempboxa >\hsize

       {#1}\ \ #2\par
   \else
       \hbox to\hsize{\hfil\box\@tempboxa\hfil}%
   \fi}
\makeatother

\newtheorem{thm}{Theorem}[section]
\newtheorem{cor}{Corollary}[section]
\newtheorem{lem}{Lemma}[section]

\newtheorem{obs}{Observation}[section]
\newtheorem{pro}{Proposition}[section]

\newcommand{\qed}{{\hfill\rule{3pt}{7pt}}}

\def\qed{\hfill \rule{4pt}{7pt}}

\begin{document}
\title{\textbf{The minimal size of graphs with given pendant-tree connectivity} \footnote{Supported by the National Science Foundation of China
(Nos. 11551001, 11161037, 11461054) and the Science Found of Qinghai
Province (No. 2014-ZJ-907).}}
\author{
\small Yaping Mao$^{1,2}$\footnote{E-mail: maoyaping@ymail.com}\\[0.2cm]
\small $^{1}$Department of Mathematics, Qinghai Normal University,\\[0.1mm]
\small $^{2}$Key Laboratory of IOT of Qinghai Province,\\[0.1mm]
\small Xining, Qinghai 810008, China}
\date{}
\maketitle

\begin{abstract}
The concept of pendant-tree $k$-connectivity $\tau_k(G)$ of a graph
$G$, introduced by Hager in 1985, is a generalization of classical
vertex-connectivity. Let $f(n,k,\ell)$ be the minimal number of
edges of a graph $G$ of order
$n$ with $\tau_k(G)=\ell \ (1\leq \ell\leq n-k)$. In this paper,
we give some exact value or sharp bounds of the parameter $f(n,k,\ell)$.\\[2mm]
\end{abstract}

\section{Introduction}

A processor network is expressed as a graph, where a node is a
processor and an edge is a communication link. Broadcasting is the
process of sending a message from the source node to all other nodes
in a network. It can be accomplished by message dissemination in
such a way that each node repeatedly receives and forwards messages.
Some of the nodes and/or links may be faulty. However, multiple
copies of messages can be disseminated through disjoint paths. We
say that the broadcasting succeeds if all the healthy nodes in the
network finally obtain the correct message from the source node
within a certain limit of time. A lot of attention has been devoted
to fault-tolerant broadcasting in networks \cite{Fragopoulou,
Hedetniemi, Jalote, Ramanathan}. In order to measure the ability of
fault-tolerance, the above path structure connecting two nodes are
generalized into some tree structures connecting more than two
nodes, see \cite{Ku, LLSun, LM1}. To show these generalizations
clearly, we must state from the connectivity in graph theory. We
divide our introduction into the following four subsections to state
the motivations and our results of this paper.

\subsection{Connectivity and $k$-connectivity}

All graphs considered in this paper are undirected, finite and
simple. We refer to the book \cite{bondy} for graph theoretical
notation and terminology not described here. For a graph $G$, let
$V(G)$, $E(G)$, $e(G)$, $\bar{G}$, $\Delta(G)$ and $\delta(G)$
denote the set of vertices, the set of edges, the size, the
complement, the maximum degree and the minimum degree of $G$,
respectively. In the sequel, let $K_{s,t}$, $K_{n}$, $W_{n}$,
$C_{n}$, and $P_n$ denote the complete bipartite graph of order
$s+t$ with part sizes $s$ and $t$, complete graph of order $n$,
wheel of order $n$, cycle of order $n$, and path of order $n$,
respectively. For any subset $X$ of $V(G)$, let $G[X]$ denote the
subgraph induced by $X$, and $E[X]$ the edge set of $G[X]$. For two
subsets $X$ and $Y$ of $V(G)$ we denote by $E_G[X,Y]$ the set of
edges of $G$ with one end in $X$ and the other end in $Y$. If
$X=\{x\}$, we simply write $E_G[x,Y]$ for $E_G[\{x\},Y]$.

Connectivity is one of the most basic concepts of graph-theoretic
subjects, both in combinatorial sense and the algorithmic sense. It
is well-known that the classical connectivity has two equivalent
definitions. The \emph{connectivity} of $G$, written $\kappa(G)$, is
the minimum order of a vertex set $S\subseteq V(G)$ such that $G-S$
is disconnected or has only one vertex. We call this definition the
`cut' version definition of connectivity. A well-known theorem of
Whitney \cite{Whitney} provides an equivalent definition of
connectivity, which can be called the `path' version definition of
connectivity. For any two distinct vertices $x$ and $y$ in $G$, the
\emph{local connectivity} $\kappa_{G}(x,y)$ is the maximum number of
internally disjoint paths connecting $x$ and $y$. Then
$\kappa(G)=\min\{\kappa_{G}(x,y)\,|\,x,y\in V(G),x\neq y\}$ is
defined to be the \emph{connectivity} of $G$. For connectivity,
Oellermann gave a survey paper on this subject; see
\cite{Oellermann2}.

Although there are many elegant and powerful results on connectivity
in graph theory, the basic notation of classical connectivity may
not be general enough to capture some computational settings. So
people want to generalize this concept. For the `cut' version
definition of connectivity, we find the above minimum vertex set
without regard the number of components of $G-S$. Two graphs with
the same connectivity may have differing degrees of vulnerability in
the sense that the deletion of a vertex cut-set of minimum
cardinality from one graph may produce a graph with considerably
more components than in the case of the other graph. For example,
the star $K_{1,n}$ and the path $P_{n+1}\ (n\geq 3)$ are both trees
of order $n+1$ and therefore connectivity $1$, but the deletion of a
cut-vertex from $K_{1,n}$ produces a graph with $n$ components while
the deletion of a cut-vertex from $P_{n+1}$ produces only two
components. Chartrand et al. \cite{Chartrand1} generalized the `cut'
version definition of connectivity. For an integer $k \ (k\geq 2)$
and a graph $G$ of order $n \ (n\geq k)$, the
\emph{$k$-connectivity} $\kappa'_k(G)$ is the smallest number of
vertices whose removal from $G$ of order $n \ (n\geq k)$ produces a
graph with at least $k$ components or a graph with fewer than $k$
vertices. Thus, for $k=2$, $\kappa'_2(G)=\kappa(G)$. For more
details about $k$-connectivity, we refer to \cite{Chartrand1, Day,
Oellermann2, Oellermann3}.

\subsection{Generalized (edge-)connectivity}

The generalized connectivity of a graph $G$, introduced by Hager
\cite{Hager}, is a natural generalization of the `path' version
definition of connectivity. For a graph $G=(V,E)$ and a set
$S\subseteq V(G)$ of at least two vertices, \emph{an $S$-Steiner
tree} or \emph{a Steiner tree connecting $S$} (or simply, \emph{an
$S$-tree}) is a such subgraph $T=(V',E')$ of $G$ that is a tree with
$S\subseteq V'$. Note that when $|S|=2$ an $S$-Steiner tree is just
a path connecting the two vertices of $S$. Two $S$-Steiner trees $T$
and $T'$ are said to be \emph{internally disjoint} if $E(T)\cap
E(T')=\varnothing$ and $V(T)\cap V(T')=S$. For $S\subseteq V(G)$ and
$|S|\geq 2$, the \emph{generalized local connectivity} $\kappa_G(S)$
is the maximum number of internally disjoint $S$-Steiner trees in
$G$, that is, we search for the maximum cardinality of edge-disjoint
trees which include $S$ and are vertex disjoint with the exception
of $S$. For an integer $k$ with $2\leq k\leq n$, \emph{generalized
$k$-connectivity} (or \emph{$k$-tree-connectivity}) is defined as
$\kappa_k(G)=\min\{\kappa_G(S)\,|\,S\subseteq V(G),|S|=k\}$, that
is, $\kappa_k(G)$ is the minimum value of $\kappa_G(S)$ when $S$
runs over all $k$-subsets of $V(G)$. Clearly, when $|S|=2$,
$\kappa_2(G)$ is nothing new but the connectivity $\kappa(G)$ of
$G$, that is, $\kappa_2(G)=\kappa(G)$, which is the reason why one
addresses $\kappa_k(G)$ as the generalized connectivity of $G$. By
convention, for a connected graph $G$ with less than $k$ vertices,
we set $\kappa_k(G)=1$. Set $\kappa_k(G)=0$ when $G$ is
disconnected. Note that the generalized $k$-connectivity and
$k$-connectivity of a graph are indeed different. Take for example,
the graph $H_1$ obtained from a triangle with vertex set
$\{v_1,v_2,v_3\}$ by adding three new vertices $u_1,u_2,u_3$ and
joining $v_i$ to $u_i$ by an edge for $1 \leq i\leq 3$. Then
$\kappa_3(H_1)=1$ but $\kappa'_3(H_1)=2$. There are many results on
the generalized connectivity, see \cite{Chartrand2, LLSun, LL, LLZ,
LM1, LM2, LM3, LM4, LMS, Okamoto}.

\subsection{Pendant-tree (edge-)connectivity}

The concept of pendant-tree connectivity \cite{Hager} was introduced
by Hager in 1985, which is specialization of generalized
connectivity (or \emph{$k$-tree-connectivity}) but a generalization
of classical connectivity. For an $S$-Steiner tree, if the degree of
each vertex in $S$ is equal to one, then this tree is called a
\emph{pendant $S$-Steiner tree}. Two pendant $S$-Steiner trees $T$
and $T'$ are said to be \emph{internally disjoint} if $E(T)\cap
E(T')=\varnothing$ and $V(T)\cap V(T')=S$. For $S\subseteq V(G)$ and
$|S|\geq 2$, the \emph{pendant-tree local connectivity} $\tau_G(S)$
is the maximum number of internally disjoint pendant $S$-Steiner
trees in $G$. For an integer $k$ with $2\leq k\leq n$,
\emph{pendant-tree $k$-connectivity} is defined as
$\tau_k(G)=\min\{\tau_G(S)\,|\,S\subseteq V(G),|S|=k\}$. When $k=2$,
$\tau_2(G)=\tau(G)$ is just the connectivity of a graph $G$. For
more details on pendant-tree connectivity, we refer to \cite{Hager,
MaoL}. Clearly, we have
$$
\left\{
\begin{array}{ll}
\tau_k(G)=\kappa_k(G),&\mbox {\rm for}~k=1,2;\\
\tau_k(G)\leq \kappa_k(G),&\mbox {\rm for}~k\geq 3.
\end{array}
\right.
$$

The relation between pendant-tree connectivity and generalized
connectivity are shown in the following Table 2.
\begin{center}
\begin{tabular}{|c|c|c|}
\hline &Pendant tree-connectivity & Generalized connectivity\\[0.1cm]
\cline{1-3}
Vertex subset & $S\subseteq V(G) \ (|S|\geq 2)$ & $S\subseteq V(G) \ (|S|\geq 2)$\\[0.1cm]
\cline{1-3} Set of Steiner trees & $\left\{
\begin{array}{ll}
\mathscr{T}_{S}=\{T_1,T_2,\cdots,T_{\ell}\}\\
S\subseteq V(T_i),\\
d_{T_i}(v)=1~for~every~v\in S\\
E(T_i)\cap E(T_j)=\varnothing,\\
\end{array}
\right.$ & $\left\{
\begin{array}{ll}
\mathscr{T}_{S}=\{T_1,T_2,\cdots,T_{\ell}\}\\
S\subseteq V(T_i),\\
E(T_i)\cap E(T_j)=\varnothing,\\
\end{array}
\right.$\\[0.05cm]
\cline{1-3}
Local parameter & $\tau(S)=\max|\mathscr{T}_{S}|$ & $\kappa(S)=\max|\mathscr{T}_{S}|$\\
\cline{1-3} Global parameter & $\tau_k(G)=\underset{S\subseteq V(G),
|S|=k}{\min} \tau(S)$ & $\kappa_k(G)=\underset{S\subseteq V(G),
|S|=k}{\min}
\kappa(S)$\\[0.05cm]
\cline{1-3}
\end{tabular}
\end{center}
\begin{center}
{Table 2. Two kinds of tree-connectivities}
\end{center}

The following two observations are easily seen.
\begin{obs}\label{obs1-1}
If $G$ is a connected graph, then $\tau_k(G)\leq \mu_k(G)\leq
\delta(G)$.
\end{obs}
\begin{obs}\label{obs1-2}
If $H$ is a spanning subgraph of $G$, then $\tau_k(H)\leq
\tau_k(G)$.
\end{obs}

In {\upshape\cite{Hager}}, Hager derived the following results.
\begin{lem}{\upshape\cite{Hager}}\label{lem1-1}
Let $G$ be a graph. If $\tau_k(G)\geq \ell$, then $\delta(G)\geq
k+\ell-1$.
\end{lem}
\begin{lem}{\upshape\cite{Hager}}\label{lem1-2}
Let $G$ be a graph. If $\tau_k(G)\geq \ell$, then $\kappa(G)\geq
k+\ell-2$.
\end{lem}
\begin{lem}{\upshape\cite{Hager}}\label{lem1-3}
Let $k,n$ be two integers with $3\leq k\leq n$, and let $K_n$ be a
complete graph of order $n$. Then
$$
\tau_k(K_n)=n-k.
$$
\end{lem}
\begin{lem}{\upshape\cite{Hager}}\label{lem1-4}
Let $K_{r,s}$ be a complete bipartite graph with $r+s$ vertices.
Then
$$
\tau_k(K_{r,s})=\max\{\min\{r-k+1, s-k+1\},0\}.
$$
\end{lem}

As a natural counterpart of the pendant-tree $k$-connectivity, we
introduced the concept of pendant-tree $k$-edge-connectivity. For
$S\subseteq V(G)$ and $|S|\geq 2$, the {\it pendant-tree local
edge-connectivity} $\mu(S)$ is the maximum number of edge-disjoint
pendant $S$-Steiner trees in $G$. For an integer $k$ with $2\leq
k\leq n$, the {\it pendant-tree $k$-edge-connectivity} $\mu_k(G)$ of
$G$ is then defined as $\mu_k(G)=\min\{\mu(S)\,|\,S\subseteq V(G) \
and \ |S|=k\}$. It is also clear that when $|S|=2$, $\mu_2(G)$ is
just the standard edge-connectivity $\lambda(G)$ of $G$, that is,
$\mu_2(G)=\lambda(G)$.

\subsection{Application background and our results}

In addition to being a natural combinatorial measure, both the
pendant-tree connectivity and the generalized connectivity can be
motivated by its interesting interpretation in practice. For
example, suppose that $G$ represents a network. If one considers to
connect a pair of vertices of $G$, then a path is used to connect
them. However, if one wants to connect a set $S$ of vertices of $G$
with $|S|\geq 3$, then a tree has to be used to connect them. This
kind of tree with minimum order for connecting a set of vertices is
usually called a Steiner tree, and popularly used in the physical
design of VLSI (see \cite{Grotschel1, Grotschel2, Sherwani}) and
computer communication networks (see \cite{Du}). Usually, one wants
to consider how tough a network can be, for the connection of a set
of vertices. Then, the number of totally independent ways to connect
them is a measure for this purpose. The generalized $k$-connectivity
can serve for measuring the capability of a network $G$ to connect
any $k$ vertices in $G$.

Let $f(n,k,\ell)$ be the minimal number of edges of a graph $G$ of
order $n$ with $\tau_k(G)=\ell \ (1\leq \ell\leq n-k)$. It is not
easy to determine the exact value of the parameter $f(n,k,\ell)$ for
a general $k  \ (3\leq k\leq n)$ and a general $\ell \ (1\leq
\ell\leq n-k)$.

In Section 2, we obtain the following result for general $k$.
\begin{thm}\label{th1-1}
Let $n,k$ be two integers with $3\leq k\leq n$ and $n\geq 15$. Then

$(1)$ $f(n,k,n-k)={n\choose{2}}$;

$(2)$ $f(n,k,n-k-1)={n\choose{2}}-2$;

$(3)$ $f(n,k,0)=n-1$;

$(4)$ $f(n,k,1)=\left\lceil \frac{kn}{2}\right \rceil$;

$(5)$ For $1\leq \ell\leq \frac{n}{2}-k+1$, we have
$$
\left\lceil \frac{1}{2}(k+\ell-1)n \right \rceil\leq f(n,k,\ell)\leq
(k+\ell-1)(n-k-\ell+1);
$$

$(6)$ For $\frac{n}{2}-k+2\leq \ell\leq n-k-2$, we have
$$
\left\lceil \frac{1}{2}(k+\ell-1)n \right \rceil\leq f(n,k,\ell)\leq
(k+\ell-1)(n-k-\ell+1)+{k+\ell-1\choose 2}.
$$

Moreover, the bounds are sharp.
\end{thm}

For $k=n,n-1,n-2,n-3,3$, we get the following results in Section
$3$.
\begin{thm}\label{th1-2}
Let $n,k$ be two integers with $3\leq k\leq n$ and $n\geq 15$. Then

$(1)$ $f(n,n,0)=n-1$;

$(2)$ $f(n,n-1,1)={n\choose{2}}$, $f(n,n-1,0)=n-1$;

$(3)$ $f(n,n-2,2)={n\choose{2}}$, $f(n,n-2,1)={n\choose{2}}-2$,
$f(n,n-2,0)=n-1$;

$(4)$ $f(n,n-3,3)={n\choose{2}}$, $f(n,n-3,2)={n\choose{2}}-2$,
$f(n,n-3,1)={n\choose{2}}-n$; $f(n,n-3,0)=n-1$.
\end{thm}

\begin{thm}\label{th1-3}
Let $n$ be an integer with $n\geq 10$. Then

$(1)$ $f(n,3,0)=n-1$;

$(2)$ $f(n,3,1)=\lceil\frac{3n}{2}\rceil$;

$(3)$ $f(n,3,2)=2n$ for $n=pq$ and $p,q\geq 3$; $2n\leq f(n,3,2)\leq
\frac{5n}{2}$ for $n=2p$; $2n\leq f(n,3,2)\leq 4n-16$ if $n$ is a
prime number.

$(4)$ $\left\lceil \frac{(\ell+2)n}{2}\right\rceil\leq
f(n,3,\ell)\leq (\ell+1)(n-\ell)+1$ for $3\leq \ell\leq
\frac{n-4}{3}$;

$(5)$ $\left\lceil \frac{(\ell+2)n}{2}\right\rceil\leq
f(n,3,\ell)\leq \left\lfloor\frac{n}{\ell+1}\right\rfloor
(\ell+1)^2+(r+1)(\ell+1)+r-1$ for $\frac{n-1}{3}\leq \ell\leq
\frac{n-r-2}{2}$ where $n\equiv r \ (mod \ \ell)$ and $r\geq 2$;
$\left\lceil \frac{(\ell+2)n}{2}\right\rceil\leq f(n,3,\ell)\leq
(\ell+2)(n-\ell-2)$ for $\frac{n-r-2}{2}\leq \ell\leq
\frac{n-4}{2}$, or $\frac{n-1}{3}\leq \ell\leq \frac{n-r-2}{2}$
where $n\equiv r \ (mod \ \ell)$ and $r=0,1$;

$(6)$ $\left\lceil \frac{(\ell+2)n}{2}\right\rceil\leq
f(n,3,\ell)\leq (\ell+2)(n-\ell-2)+{\ell+2\choose 2}$ for
$\frac{n-2}{2}\leq \ell\leq n-6$;

$(7)$ $f(n,3,n-5)={n\choose 2}-\lfloor\frac{n-4}{2}\rfloor$;

$(8)$ $f(n,3,n-4)={n\choose 2}-2$;

$(9)$ $f(n,3,n-3)={n\choose 2}$.
\end{thm}

\section{For general $k$ and $\ell$}

Given a vertex $x$ and a set $U$ of vertices, an \emph{$(x,U)$-fan}
is a set of paths from $x$ to $U$ such that any two of them share
only the vertex $x$. The size of an $(x,U)$-fan is the number of
internally disjoint paths from $x$ to $U$.
\begin{lem}{\upshape (Fan Lemma, \cite{West}, p-170)}\label{lem2-1}
A graph is $k$-connected if and only if it has at least $k+1$
vertices and, for every choice of $x$, $U$ with $|U|\geq k$, it has
an $(x,U)$-fan of size $k$.
\end{lem}

\begin{cor}\label{cor2-1}
Let $G$ be a graph with $\kappa(G)=k$. Then $\tau_k(G)\geq 1$.
\end{cor}
\begin{pf}
Since $\kappa(G)=k$, it follows from Lemma \ref{lem2-1} that for any
$S\subseteq V(G)$ and $|S|=k$, there exist an $(x,S)$-fan, where
$x\in V(G)-S$. So there is a pendant $S$-Steiner tree in $G$, and
hence $\tau(S)\geq 1$. From the arbitrariness of $S$, we have
$\tau_k(G)\geq 1$, as desired.\qed
\end{pf}

In \cite{MaoLai}, Mao and Lai obtained the following upper bound of
$\tau_k(G)$.
\begin{lem}{\upshape \cite{MaoLai}}\label{lem2-2}
Let $G$ be a connected graph of order $n\geq 6$. Let $k$ be an
integer with $3\leq k\leq n$. Then
$$
\tau_k(G)\leq \delta(G)-k+1.
$$
Moreover, the upper bound is sharp.
\end{lem}

The Harary graph $H_{n,d}$ is constructed by arranging the $n$
vertices in a circular order and spreading the $d$ edges around the
boundary in a nice way, keeping the chords as short as possible.
{\it Harary graph}\index{Harary graph} $H_{n,d}$ is a $d$-connected
graph on $n$ vertices, and the structure of $H_{n,d}$ depends on the
parities of $d$ and $n$; see \cite{West}.
\begin{itemize}
\item[] \textbf{Case 1:} $d$ even. Let $d=2r$. Then $H_{n,2r}$ is constructed as
follows. It has vertices $0,1,\cdots,n-1$ and two vertices $i$ and
$j$ are jointed if $i-r\leq j\leq i+r$ (where addition is taken
modulo $n$).

\item[] \textbf{Case 2:} $d$ odd, $n$ even. Let $d=2r+1$. Then $H_{n,2r+1}$ is
constructed by first drawing $H_{n,2r}$ and then adding edges
joining vertex $i$ to vertex $i+\frac{n}{2}$ for $1\leq i\leq
\frac{n}{2}$.

\item[] \textbf{Case 3:} $d$ odd, $n$ even. Let $d=2r+1$. Then $H_{n,2r+1}$ is
constructed by first drawing $H_{n,2r}$ and then adding edges
joining vertex $0$ to vertices $\frac{n-1}{2}$ and $\frac{n+1}{2}$
and $i$ to vertex $i+\frac{n+1}{2}$ for $1\leq i\leq \frac{n-1}{2}$.
\end{itemize}

\begin{lem}{\upshape \cite{West}}\label{lem2-3}
Let $H_{n,k}$ be the Harary graph of order $n\geq 6$. Then
$$
\kappa(H_{n,k})=k, \ e(H_{n,k})=\left\lceil\frac{kn}{2}\right\rceil.
$$
\end{lem}

For general $k$ and $\ell$, we give a sharp lower bound of
$f(n,k,\ell)$.
\begin{pro}\label{pro2-1}
Let $n,k$ be two integers with $3\leq k\leq n$ and $n\geq 15$. Then
$$
f(n,k,\ell)\geq \left\lceil \frac{1}{2}(k+\ell-1)n \right \rceil
$$
for $1\leq \ell\leq n-k-2$. Moreover, the bound is sharp.
\end{pro}
\begin{pf}
Since $\tau_k(G)=\ell \ (1\leq \ell\leq n-k-2)$, it follows from
Lemma \ref{lem2-2} that $\delta(G)\geq k+\ell-1$. Denote by $X$ the
set of vertices of degree $k+\ell-1$ in $G$. Set $Y=V(G)\setminus
X$. Let $m',m''$ be the number of edges of $G[X],G[Y]$,
respectively. Then
$$
|E_G[X,Y]|=(k+\ell-1)|X|-2m' \eqno (1)
$$
and
\begin{eqnarray*}
~~~~~~~~~~~~~~~~~~~~e(G)&=&e(G[X])+e(G[Y])+|E_G[X,Y]|\\
~~~~~~~~~~~~~~~~~~~~&=&m'+m''+(k+\ell-1)|X|-2m'\\
~~~~~~~~~~~~~~~~~~~~&=&(k+\ell-1)|X|+(m''-m').~~~~~~~~~~~~~~~~~~~~~~~~~~~~~~~~~~~~~~~~~~(2)
\end{eqnarray*}
Since every vertex in $Y$ has degree at least $k+\ell$ in $G$, it
follow that
\begin{eqnarray*}
\sum_{v\in Y}d(v)&=&2m''+|E_G[X,Y]|\\
&=&2m''+(k+\ell-1)|X|-2m'\\
&\geq&(k+\ell)|Y|,
\end{eqnarray*}
and hence
$$
m''-m'\geq \frac{1}{2}(k+\ell-1)(|Y|-|X|)+\frac{1}{2}|Y|.\eqno (3)
$$
From $(1),(2),(3)$, we have
$$
e(G)=\frac{1}{2}(k+\ell-1)(|Y|+|X|)+\frac{1}{2}|Y|\geq
\frac{1}{2}(k+\ell-1)n.
$$
Since the number of edges is an integer, it follows that
$$
e(G)\geq
\left\lceil \frac{1}{2}(k+\ell-1)n \right\rceil,
$$
as desired.\qed
\end{pf}

To show the sharpness of the lower bound, we consider the following
example.

\noindent \textbf{Example 1.} Let $H_{n,k}$ be the Harary graph of
order $n\geq 6$. Then $\kappa(H_{n,k})=\delta(H_{n,k})=k$. From
Lemma \ref{lem2-2}, we have $\tau_k(H_{n,k})\leq
\delta(H_{n,k})-k+1=1$. From Corollary \ref{cor2-1}, we have
$\tau_k(H_{n,k})\geq 1$, and hence $\tau_k(H_{n,k})=1$. From Lemma
\ref{lem2-3},
$e(H_{n,k})=\left\lceil\frac{kn}{2}\right\rceil=f(n,k,1)$. So the
lower bound is sharp for $\ell=1$.

The following corollary is immediate from Proposition \ref{pro2-1}
and Example 1.
\begin{cor}\label{cor2-2}
Let $n,k$ be two integers with $3\leq k\leq n$ and $n\geq 5$. Then
$$
f(n,k,1)=\left\lceil \frac{kn}{2}\right \rceil.
$$
\end{cor}

For general $k$ and $\ell \ (1\leq \ell\leq \frac{n}{2}-k+1)$, we
give a sharp upper bound of $f(n,k,\ell)$.
\begin{pro}\label{pro2-2}
Let $n,k,\ell$ be three integers with $3\leq k\leq n$ and $1\leq
\ell\leq \frac{n}{2}-k+1$. Then
$$
f(n,k,\ell)\leq (k+\ell-1)(n-k-\ell+1).
$$
Moreover, the bound is sharp.
\end{pro}
\begin{pf}
Let $G=K_{k+\ell-1,n-k-\ell+1}$. Since $1\leq \ell\leq
\frac{n}{2}-k+1$, it follows that $n-k-\ell+1\geq k+\ell-1$. From
Lemma \ref{lem1-4}, we have $\tau_k(G)=\ell$. So $f(n,k,\ell)\leq
(k+\ell-1)(n-k-\ell+1)$, as desired. \qed
\end{pf}

The {\it join} or {\it complete product} $G\vee H$ of two disjoint
graphs $G$ and $H$ is the graph with vertex set $V(G)\cup V(H)$ and
edge set $E(G)\cup E(H)\cup \{uv\,|\, u\in V(G), v\in V(H)\}$. For
classical connectivity, the following result is well-known.
\begin{lem}{\upshape \cite{bondy}}\label{lem2-4}
Let $G$ and $H$ be two graphs. Then
$$
\kappa(G\vee H)=\min\{|V(G)|+\kappa(H),|V(H)|+\kappa(G)\}.
$$
\end{lem}

For general $k$ and $\ell \ (\frac{n}{2}-k+2\leq \ell\leq n-k)$, we
can also give a sharp upper bound of $f(n,k,\ell)$.
\begin{pro}\label{pro2-3}
Let $n,k,\ell$ be three integers with $3\leq k\leq n$ and
$\frac{n}{2}-k+2\leq \ell\leq n-k$. Then
$$
f(n,k,\ell)\leq (k+\ell-1)(n-k-\ell+1)+{k+\ell-1\choose 2}.
$$
Moreover, the bound is sharp.
\end{pro}
\begin{pf}
Let $G=K_{k+\ell-1}\vee (n-k-\ell+1)K_1$. Since
$\delta(G)=k+\ell-1$, it follows from Lemma \ref{lem2-2} that
$\tau_k(G)\leq \delta(G)-k+1=\ell$. It suffices to show that
$\tau_k(G)\geq \ell$. Let
$V(K_{k+\ell-1})=\{y_1,y_2,\ldots,y_{k+\ell-1}\}$ and
$X=V(G)-V(K_{k+\ell-1})=\{x_1,x_2,\ldots,x_{n-k-\ell+1}\}$. For any
$S\subseteq V(G)$ and $|S|=k$, if $S\subseteq X$, then the trees
induced by the edges in $E_G[y_i,S] \ (1\leq i\leq k+\ell-1)$ are
$k+\ell-1$ internally disjoint pendant $S$-Steiner trees, and hence
$\tau(S)\geq k+\ell-1$. If $S\subseteq Y$, then there are $\ell-1$
internally disjoint pendant $S$-Steiner trees in $G[Y]$. These trees
and the trees induced by the edges in $E_G[x_i,S] \ (1\leq i\leq
n-k-\ell+1)$ are $n-k$ internally disjoint pendant $S$-Steiner
trees, and hence $\tau(S)\geq n-k\geq \ell$. Suppose $S\cap X\neq
\varnothing$ and $S\cap V(K_{k+\ell-1})\neq \varnothing$. Without
loss of generality, let $S\cap X=\{x_1,x_2,\ldots,x_{r}\}$ and
$S\cap V(K_{k+\ell-1})=\{y_1,y_2,\ldots,y_{k-r}\}$. Then the trees
induced by the edges in $\{x_iy_j\,|\,1\leq i\leq r, \ k-r+1\leq
j\leq k+\ell-1\}\cup \{y_iy_j\,|\,1\leq i\leq k-r, \ k-r+1\leq j\leq
k+\ell-1\}$ are $\ell-1+r\geq \ell$ internally disjoint pendant
$S$-Steiner trees, and hence $\tau(S)\geq n-1$. We conclude that
$\tau(S)\geq n-1$ for any $S\subseteq V(G)$ and $|S|=k$. From the
arbitrariness of $S$, we have $\tau_k(G)\geq \ell$, and hence
$\tau_k(G)=\ell$. So $f(n,k,\ell)\leq
(k+\ell-1)(n-k-\ell+1)+{k+\ell-1\choose 2}$. \qed
\end{pf}

In \cite{MaoL}, Mao and Lai characterized graphs with
$\tau_k(G)=n-k,n-k-1,n-k-2$, respectively.
\begin{lem}{\upshape \cite{MaoL}}\label{lem2-5}
Let $k,n$ be two integers with $3\leq k\leq n$ and $n\geq 4$, and
let $G$ be a connected graph. Then $\tau_k(G)=n-k$ if and only if
$G$ is a complete graph of order $n$.
\end{lem}
\begin{lem}{\upshape \cite{MaoL}}\label{lem2-6}
Let $k,n$ be two integers with $3\leq k\leq n$ and $n\geq 7$, and
let $G$ be a connected graph. Then $\tau_k(G)=n-k-1$ if and only if
$\bar{G}=r K_2\cup (n-2r)K_1\ (r=1,2)$.
\end{lem}

From $\ell=n-k,n-k-1,0$, we can derive the exact value of
$f(n,k,\ell)$.
\begin{pro}\label{pro2-4}
Let $n,k$ be two integers with $3\leq k\leq n$ and $n\geq 15$. Then

$(i)$ $f(n,k,n-k)={n\choose{2}}$;

$(ii)$ $f(n,k,n-k-1)={n\choose{2}}-2$;

$(iii)$ $f(n,k,0)=n-1$.
\end{pro}
\begin{pf}
$(i)$ From Lemma \ref{lem2-5}, $\tau_k(G)=n-k$ if and only if $G$ is
a complete graph of order $n$. So $f(n,k,n-k)={n\choose{2}}$.

$(ii)$ From Lemma \ref{lem2-6}, $\tau_k(G)=n-k-1$ if and only if
$\bar{G}=r K_2\cup (n-2r)K_1\ (r=1,2)$. Then
$f(n,k,n-k-1)={n\choose{2}}-2$.

$(iii)$ The tree $T_n$ with $n$ vertices is the graph such that
$\tau_k(T_n)=0$ with the minimal number of edges. So
$f(n,k,1)=n-1$.\qed
\end{pf}

The results in Theorem \ref{th1-1} follow from Propositions
\ref{pro2-1}, \ref{pro2-2}, \ref{pro2-3} and \ref{pro2-4}.

\section{For fixed $k$ and general $\ell$}

The {\it Cartesian product} $G\Box H$ of two graphs $G$ and $H$, is
the graph with vertex set $V(G)\times V(H)$, in which two vertices
$(u,v)$ and $(u',v')$ are adjacent if and only if $u=u'$ and
$(v,v')\in E(H)$, or $v=v'$ and $(u,u')\in E(G)$. Clearly, $|E(G\Box
H)|=|E(H)||V(G)|+|E(G)||V(H)|$.

The {\it lexicographic product} $G\circ H$ of graphs $G$ and $H$ has
the vertex set $V(G\circ H)=V(G)\times V(H)$, and two vertices
$(u,v),(u',v')$ are adjacent if $uu'\in E(G)$, or if $u=u'$ and
$vv'\in E(H)$. It is easy to see that $|E(G\circ
H)|=|E(H)||V(G)|+|E(G)||V(H)|^2$.

In this section, we study $f(n,k,\ell)$ for $k=n,n-1,n-2,n-3,3$  and
general $\ell$.

\subsection{For $n-3\leq k\leq n$ and general $\ell$}

Mao and Lai obtained the following results.
\begin{lem}{\upshape \cite{MaoL}}\label{lem3-1}
Let $G$ be a graph of order $n$. Then $\tau_n(G)=0$ if and only if
$G$ is a graph of order $n$.
\end{lem}
\begin{lem}{\upshape \cite{MaoL}}\label{lem3-2}
Let $G$ be a connected graph of order $n$. Then

$(1)$ $\tau_{n-1}(G)=1$ if and only if $G$ is a complete graph of
order $n$.

$(2)$ $\tau_{n-1}(G)=0$ if and only if $G$ is not a complete graph
of order $n$.
\end{lem}
\begin{lem}{\upshape \cite{MaoL}}\label{lem3-3}
Let $G$ be a connected graph of order $n$. Then

$(1)$ $\tau_{n-2}(G)=2$ if and only if $G$ is a complete graph of
order $n$.

$(2)$ $\tau_{n-2}(G)=1$ if and only if $G=K_n\setminus M$ and $1\leq
|M|\leq 2$, where $M$ is a matching of $K_n$ for $n\geq 7$.

$(3)$ $\tau_{n-2}(G)=0$ if and only if $G$ is one of the other
graphs.
\end{lem}

Graphs with $\tau_{n-3}(G)=\ell \ (0\leq \ell\leq 3)$ can be also
characterized.
\begin{pro}\label{pro3-1}
Let $G$ be a connected graph of order $n\geq 9$. Then

$(1)$ $\tau_{n-3}(G)=3$ if and only if $G$ is a complete graph of
order $n$.

$(2)$ $\tau_{n-3}(G)=2$ if and only if $G=K_n\setminus M$ and $1\leq
|M|\leq 2$, where $M$ is a matching of $K_n$.

$(3)$ $\tau_{n-3}(G)=1$ if and only if $1\leq \Delta(\bar{G})\leq
2$, and $|M|\geq 3$ if $M$ is a matching of $\bar{G}$.

$(4)$ $\tau_{n-3}(G)=0$ if and only if $G$ is one of the other
graphs.
\end{pro}
\begin{pf}
From Lemmas \ref{lem2-5} and \ref{lem2-6}, $(1)$ and $(2)$ is true.
We only need to show that $(3)$ is true. Suppose that
$\tau_{n-3}(G)=1$. We claim that $1\leq \Delta(\bar{G})\leq 2$.
Assume, to the contrary, that $\Delta(\bar{G})\geq 3$. Then there
exist four vertices $u,v,w,u_1$ such that $u_1u,u_1v,u_1w\notin
E(G)$. We choose $S\subseteq V(G)$ and $|S|=n-3$ such that
$V(G)-S=\{u,v,w\}$. Then $u_1\in S$. Observe that there is no
pendant $S$-Steiner tree in $G$. Therefore, $\tau_{n-3}(G)=0$, a
contradiction. So $1\leq \Delta(\bar{G})\leq 2$. From $(2)$, if $M$
is a matching of $\bar{G}$, then $|M|\geq 3$.

Conversely, we suppose $1\leq \Delta(\bar{G})\leq 2$, and $|M|\geq
3$ if $M$ is a matching of $\bar{G}$. For any $S\subseteq V(G)$ and
$|S|=n-3$, there exist three vertices, say $u,v,w$, such that
$\bar{S}=V(G)-S=\{u,v,w\}$. Since $1\leq \Delta(\bar{G})\leq 2$, it
follows that $1\leq d_{\bar{G}}(u),d_{\bar{G}}(v),d_{\bar{G}}(w)\leq
2$. Set $S=\{u_1,u_2,\ldots,u_{n-3}\}$. If
$d_{\bar{G}[\bar{S}]}(u)=2$, then $uu_i\in E(G) \ (1\leq i\leq
n-3)$, and hence the tree induced by the edges in $E_G[u,S]$ is a
pendant $S$-Steiner tree, and hence $\tau(S)\geq 1$. The same is
true when $d_{\bar{G}[\bar{S}]}(v)=2$ or
$d_{\bar{G}[\bar{S}]}(w)=2$. From now on, we suppose that
$d_{\bar{G}[\bar{S}]}(u)\leq 1$, $d_{\bar{G}[\bar{S}]}(v)\leq 1$ and
$d_{\bar{G}[\bar{S}]}(w)\leq 1$. Clearly, we have $G[\bar{S}]=P_3$
or $G[\bar{S}]=K_3$. Suppose $G[\bar{S}]=P_3$. Without loss of
generality, let $uv,uw\in E(G)$. Since $1\leq \Delta(\bar{G})\leq
2$, there are at most two vertices in $S$, say $u_1,u_2$, such that
$uu_1,uu_2\notin E(G)$. Since $1\leq \Delta(\bar{G})\leq 2$, it
follows that $u_iv\in E(G)$ or $u_iw\in E(G)$ for $i=1,2$. If
$u_1w,u_2w\in E(G)$, then the tree induced by the edges in
$\{uu_i\,|\,3\leq i\leq n-3\}\cup \{u_1w,u_2w\}$ is a pendant
$S$-Steiner tree, and hence $\tau(S)\geq 1$. If $u_1v,u_2w\in E(G)$,
then the tree induced by the edges in $\{uu_i\,|\,3\leq i\leq
n-3\}\cup \{u_1v,u_2w,uv,vw\}$ is a pendant $S$-Steiner tree, and
hence $\tau(S)\geq 1$. From the arbitrariness of $S$, we have
$\tau_{n-3}(G)\geq 1$. Since $|M|\geq 3$ if $M$ is a matching of
$\bar{G}$, it follows that $\tau_{n-3}(G)\leq 1$. So
$\tau_{n-3}(G)=1$.\qed
\end{pf}

Theorem \ref{th1-2} is immediate from Lemmas \ref{lem3-1},
\ref{lem3-2}, \ref{lem3-3} and Proposition \ref{pro3-1}.

\subsection{For $k=3$ and general $\ell$}

In \cite{Spacapan}, \u{S}pacapan obtained the following result.
\begin{lem}{\upshape\cite{Spacapan}}\label{lem3-4}
Let $G$ and $H$ be two nontrivial graphs. Then
$$
\kappa(G\Box
H)=\min\{\kappa(G)|V(H)|,\kappa(H)|V(G)|,\delta(G)+\delta(H)\}.
$$
\end{lem}

For $k=3$ and $\ell=2$, we have the following.
\begin{pro}\label{pro3-2}
Let $n$ be an integer with $n\geq 10$.

$(1)$ If $n=pq$ and $p,q\geq 3$, then $f(n,3,2)=2n$;

$(2)$ If $n=2p$, then $2n\leq f(n,3,2)\leq \frac{5n-8}{2}$;

$(3)$ If $n$ is a prime number, then $2n\leq f(n,3,2)\leq 4n-16$.
\end{pro}
\begin{pf}
$(1)$ If $n=pq$ and $p,q\geq 3$, then we consider the graph
$G=C_p\Box C_q$. From Lemma \ref{lem3-4}, $\kappa(G)=\delta(G)=4$.
From Lemma \ref{lem2-2}, we have $\tau_3(G)\leq \delta(G)-2=2$. We
only need to show that $\tau_3(G)\geq 2$. Let
$C_{p,1},C_{p,2},\ldots, C_{p,q}$ be cycles in $G$ corresponding to
$C_{p}$, and let $C_{1,q}',C_{2,q}',\ldots, C_{p,q}'$ be cycles in
$G$ corresponding to $C_{q}$. For any $S\subseteq V(G)$ and $|S|=3$,
it suffices to show $\tau(S)\geq 2$. Let $S=\{x,y,z\}$.

If there exists some cycle $C_{p,i}$ such that $|S\cap
V(C_{p,i})|=3$, then $x,y,z\in V(C_{p,i})$. Let $x',y',z'$ be the
vertices in $C_{p,i-1}$ corresponding to $x,y,z$ in $C_{p,i}$, and
let $x'',y'',z''$ be the vertices in $C_{p,i+1}$ corresponding to
$x,y,z$ in $C_{p,i}$, respectively. Then the subgraph induced by the
edges in $E(C_{p,i-1})\cup \{xx',yy',zz'\}$ contains a pendant
$S$-Steiner tree, and the subgraph induced by the edges in
$E(C_{p,i+1})\cup \{xx'',yy'',zz''\}$ contains a pendant $S$-Steiner
tree. Note that these trees are internally disjoint. Then
$\tau(S)\geq 2$.

Suppose that there exists some cycle $C_{p,i}$ such that $|S\cap
V(C_{p,i})|=2$. Then there exists another cycle $C_{p,j}$ such that
$|S\cap V(C_{p,j})|=1$. Without loss of generality, let $|S\cap
V(C_{p,1})|=2$ and $|S\cap V(C_{p,2})|=1$. Then $x,y\in V(C_{p,1})$
and $z\in V(C_{p,2})$. Let $x',y'$ be the vertices corresponding to
$x,y$ in $C_{p,1}$, and let $z'$ be the vertices corresponding to
$z$ in $C_{p,2}$. Suppose $z'\notin \{x,y\}$. Then the vertices
$x,y,z'$ divide the cycle $C_{p,1}$ into three paths, say
$P_1,P_2,P_3$. Similarly, the vertices $x',y',z$ divide the cycle
$C_{p,2}$ into three paths, say $Q_1,Q_2,Q_3$. Then the tree induced
by the edges in $E(P_{2})\cup E(P_{3})\cup \{zz'\}$ and the tree
induced by the edges in $E(Q_{1})\cup E(Q_{3})\cup \{xx',yy'\}$ are
two internally disjoint pendant $S$-Steiner trees in $G$, and hence
$\tau(S)\geq 2$. Suppose $z'\in \{x,y\}$. Without loss of
generality, let $y'=z$. Since $C_{p,1}$ is a cycle, it follows that
there exists a vertex $w\notin \{x,y\}$. Let $x',w'$ be the vertices
in $C_{p,2}$ corresponding to $x,w$ in $C_{p,1}$, and let
$x'',y'',w''$ be the vertices in $C_{p,3}$ corresponding to $x,y,w$
in $C_{p,1}$, respectively. Then the vertices $x,y,w$ divide the
cycle $C_{p,1}$ into three paths, say $P_1,P_2,P_3$. Similarly, the
vertices $x',z,w'$ divide the cycle $C_{p,2}$ into three paths, say
$Q_1,Q_2,Q_3$, and the vertices $x'',y'',w''$ divide the cycle
$C_{p,3}$ into three paths, say $R_1,R_2,R_3$. Then the tree induced
by the edges in $E(P_{2})\cup E(P_{3})\cup E(Q_{2})\cup \{ww'\}$ and
the tree induced by the edges in $E(Q_{1})\cup E(R_{1})\cup E(R)\cup
\{xx',x'x''\}$ are two internally disjoint pendant $S$-Steiner trees
in $G$, where $R$ is a path in some $C_{j,q}$ connecting $y$ and
$y''$, and hence $\tau(S)\geq 2$.

Suppose that there exist three cycles $C_{p,i},C_{p,j},C_{p,k}$ such
that $|S\cap V(C_{p,i})|=|S\cap V(C_{p,j})|=|S\cap V(C_{p,k})|=1$.
Without loss of generality, let $|S\cap V(C_{p,1})|=|S\cap
V(C_{p,2})|=|S\cap V(C_{p,3})|=1$. Let $y',z'$ be the vertices
corresponding to $y,z$ in $C_{p,1}$, $x',z''$ be the vertices
corresponding to $x,z$ in $C_{p,2}$ and $x'',y''$ be the vertices
corresponding to $x,y$ in $C_{p,3}$. If $x,y',z'$ are distinct
vertices in $C_{p,1}$, then the vertices $x,y',z'$ divide the cycle
$C_{p,1}$ into three paths $P_1,P_2,P_3$, and the vertices
$x',y,z''$ divide the cycle $C_{p,2}$ into three paths
$Q_1,Q_2,Q_3$, and the vertices $x'',y'',z$ divide the cycle
$C_{p,3}$ into three paths $R_1,R_2,R_3$. Then the tree induced by
the edges in $E(P_{1})\cup E(P_{2})\cup \{yy',z'z'',zz''\}$ and the
tree induced by the edges in $E(Q_{1})\cup E(R_{3})\cup
\{xx',x'x''\}$ are two internally disjoint pendant $S$-Steiner trees
in $G$, and hence $\tau(S)\geq 2$. Suppose that two of $x, y',z'$
are the same vertex in $C_{p,1}$. Without loss of generality, let
$x=y'$. Since $C_{p,1}$ is a cycle, it follows that there exists a
vertex $v$ such that $v\notin \{x,z'\}$. Let $v',v''$ be the
vertices in $C_{p,2},C_{p,3}$ corresponding to $v$, respectively.
Then the vertices $x,z',v$ divide the cycle $C_{p,1}$ into three
paths $P_1,P_2,P_3$, and the vertices $y,z'',v'$ divide the cycle
$C_{p,2}$ into three paths $Q_1,Q_2,Q_3$, and the vertices
$x'',z,v''$ divide the cycle $C_{p,3}$ into three paths
$R_1,R_2,R_3$. Then the tree induced by the edges in $E(P_{3})\cup
E(Q_{3})\cup E(R_{2})\cup \{vv',v'v''\}$ and the tree induced by the
edges in $E(P_{1})\cup E(Q_{1})\cup \{z'z'',zz''\}$ are two
internally disjoint pendant $S$-Steiner trees in $G$, and hence
$\tau(S)\geq 2$. Suppose that $x,y',z'$ are the same vertex in
$C_{p,1}$. Since $C_{p,1}$ is a cycle, it follows that there exists
two verties $u,v$ such that $u\neq x$ and $v\neq x$. Let $u',v'$ be
the vertices in $C_{p,2}$ corresponding to $u,v$, respectively. Let
$u'',v''$ be the vertices in $C_{p,3}$ corresponding to $u,v$,
respectively. Then the vertices $x,u,v$ divide the cycle $C_{p,1}$
into three paths $P_1,P_2,P_3$, and the vertices $y,u',v'$ divide
the cycle $C_{p,2}$ into three paths $Q_1,Q_2,Q_3$, and the vertices
$z,u'',v''$ divide the cycle $C_{p,3}$ into three paths
$R_1,R_2,R_3$. Then the tree induced by the edges in $E(P_{1})\cup
E(Q_{1})\cup E(R_{1})\cup \{uu',u'u''\}$ and the tree induced by the
edges in $E(P_{3})\cup E(Q_{3})\cup E(R_{3})\cup \{vv',v'v''\}$ are
two internally disjoint pendant $S$-Steiner trees in $G$, and hence
$\tau(S)\geq 2$.

From the argument, we conclude that $\tau_3(G)=2$, and hence
$f(n,3,2)\leq 2pq=2n$. From Proposition \ref{pro2-1}, $f(n,3,2)\geq
2n$. So $f(n,3,2)=2n$.

$(2)$ Let $G=W_p\Box P_2$. From Lemma \ref{lem3-4}, we have
$\delta(G)=\kappa(G)=4$. From Lemma \ref{lem2-2}, $\tau_3(G)\leq
\delta(G)-k+1=2$. We only need to show that $\tau_3(G)\geq 2$. It
suffices to prove that $\tau(S)\geq 2$ for any $S\subseteq V(G)$ and
$|S|=3$. Let $S=\{x,y,z\}$. Let $W_{p,1},W_{p,2}$ be the two wheels
in $G$ corresponding to $W_{p}$. Suppose $S\subseteq V(W_{p,i})$
where $i=1,2$. Without loss of generality, let $S\subseteq
V(W_{p,1})$. Let $x',y',z'$ be the vertices in $W_{p,2}$
corresponding to $x,y,z$ in $W_{p,1}$, respectively. Since
$\kappa(W_{p,1})=3$, it follows from Corollary \ref{cor2-1} that
$\tau_3(W_{p,1})\geq 1$, and hence there exists a pendant
$S$-Steiner tree in $W_{p,1}$, say $T$. Let $T'$ be the tree in
$W_{p,2}$ corresponding to $T$ in $W_{p,1}$. Then the tree induced
by the edges in $E(T')\cup \{xx',yy',zz'\}$ and the tree $T$ are two
pendant $S$-Steiner trees in $G$. Then $\tau(S)\geq 2$. Suppose
$|S\cap V(W_{p,1})|=2$ or $|S\cap V(W_{p,2})|=2$. Without loss of
generality, let $|S\cap V(W_{p,1})|=2$. Then $|S\cap V(W_{p,2})|=1$.
Let $x',y'$ be the vertices in $W_{p,2}$ corresponding to $x,y$ in
$W_{p,1}$, and let $z'$ be the vertex in $W_{p,1}$ corresponding to
$z$ in $W_{p,2}$. Let $w,w'$ be the centers of the wheels
$W_{p,1},W_{p,2}$, respectively. Then the tree $T_1$ induced by the
edges in $\{xw,yw,ww',w'z\}$ is a pendant $S$-Steiner tree, and the
tree $T_2$ induced by the edges in $\{xx',yy'\}$ is a pendant
$S$-Steiner tree. Since $T_1$ and $T_2$ are disjoint, it follows
that $\tau(S)\geq 2$. From the arbitrariness of $S$, we have
$\tau_3(G)\geq 2$, and hence $\tau_3(G)=2$. So $f(n,3,2)\leq
5p-4=\frac{5n-8}{2}$. From Proposition \ref{pro2-1}, $f(n,3,2)\geq
2n$.

$(3)$ From $(2)$ of Theorem \ref{th1-1}, if $n$ is a prime number,
then $2n\leq f(n,3,2)\leq 4n-16$.\qed
\end{pf}

For $k=3$ and general $\ell \ (3\leq \ell \leq \frac{n-4}{3})$, we
give the upper and lower bound of $f(n,3,\ell)$.
\begin{pro}\label{pro3-3}
Let $n,\ell$ be two integers with $3\leq \ell \leq \frac{n-4}{3}$
and $n\geq 15$. Then
$$
\left\lceil \frac{(\ell+2)n}{2}\right\rceil\leq f(n,3,\ell)\leq
(\ell+1)(n-\ell)+1.
$$
\end{pro}
\begin{pf}
Let $H=\ell K_{1}\vee C_{n-\ell}$. Set
$V(C_{n-\ell})=\{y_1,y_2,\ldots,y_{n-\ell}\}$. Let $G$ be a graph
obtained from $H$ by adding an edge
$y_1y_{\lfloor\frac{n-\ell}{2}\rfloor}$. Then $\delta(G)=\ell+2$.
From Lemma \ref{lem2-2}, $\tau_3(G)\leq \delta(G)-2=\ell$. We only
need to show $\tau_3(G)\geq \ell$. Set
$X=V(G)-V(C_{n-\ell})=\{x_1,x_2,\ldots,x_{\ell}\}$. For any
$S\subseteq V(G)$ and $|S|=3$, it suffices to show $\tau(S)\geq
\ell$. Suppose $S\subseteq V(C_{n-\ell})$. Without loss of
generality, let $S=\{y_1,y_2,y_3\}$. Then the trees induced by the
edges in $\{x_iy_1,x_iy_2,x_iy_3\} \ (1\leq i\leq \ell)$ are $\ell$
internally disjoint pendant $S$-Steiner trees in $G$, and hence
$\tau(S)\geq \ell$. Suppose $S\subseteq X$. Without loss of
generality, let $S=\{x_1,x_2,x_3\}$. Then the trees induced by the
edges in $\{y_ix_1,y_ix_2,y_ix_3\} \ (1\leq i\leq n-\ell)$ are
internally disjoint $n-\ell\geq 2\ell+4$ pendant $S$-Steiner trees
in $G$, and hence $\tau(S)\geq \ell$. Suppose $|S\cap X|=2$. Then
$|S\cap V(C_{n-\ell})|=1$. Without loss of generality, let
$S=\{x_1,x_2,y_1\}$. Then the trees induced by the edges in
$\{y_ix_1,y_ix_2,x_iy_i,x_iy_1\} \ (3\leq i\leq \ell)$, the tree
induced by the edges in $\{y_2x_1,y_2x_2,y_2y_1\}$ and the tree
induced by the edges in
$\{y_{n-\ell}x_1,y_{n-\ell}x_2,y_{n-\ell}y_1\}$ are $\ell$
internally disjoint pendant $S$-Steiner trees in $G$, and hence
$\tau(S)\geq \ell$. Suppose $|S\cap X|=1$. Then $|S\cap
V(C_{n-\ell})|=2$. Without loss of generality, let
$S=\{x_1,y_i,y_j\}$. If $1\leq i<j\leq
\lfloor\frac{n-\ell}{2}\rfloor$, then the trees induced by the edges
in
$\{x_ty_{\lfloor\frac{n-\ell}{2}\rfloor+t-1},y_ix_t,y_jx_t,x_1y_{\lfloor\frac{n-\ell}{2}\rfloor+t-1}\}
\ (2\leq t\leq \ell)$ and the tree induced by the edges in
$\{y_1x_{\lfloor\frac{n-\ell}{2}\rfloor},x_1y_1\}\cup
\{y_ry_{r+1}\,|\,1\leq r\leq i-1\}\cup \{y_sy_{s+1}\,|\,j\leq s\leq
\lfloor\frac{n-\ell}{2}\rfloor-1\}$ are $\ell$ internally disjoint
pendant $S$-Steiner trees in $G$, and hence $\tau(S)\geq \ell$. The
same is true for $\lfloor\frac{n-\ell}{2}\rfloor+1\leq i<j\leq
n-\ell$. Suppose $1\leq i\leq \lfloor\frac{n-\ell}{2}\rfloor$ and
$\lfloor\frac{n-\ell}{2}\rfloor+1\leq j\leq n-\ell$. Without loss of
generality, let $|j-i|\geq \lceil\frac{n-\ell}{2}\rceil$. Then the
trees induced by the edges in
$\{x_ty_{i+t},y_ix_t,y_jx_t,x_1y_{i+t}\} \ (2\leq t\leq \ell)$ and
the tree induced by the edges in $\{x_1y_1\}\cup
\{y_ry_{r+1}\,|\,1\leq r\leq i-1\}\cup \{y_sy_{s+1}\,|\,j\leq s\leq
\lfloor\frac{n-\ell}{2}\rfloor-1\}$ are $\ell$ internally disjoint
pendant $S$-Steiner trees in $G$, and hence $\tau(S)\geq \ell$. From
the argument, we conclude that $\tau(S)\geq \ell$ for any
$S\subseteq V(G)$ and $|S|=3$. So $\tau_3(G)=\ell$ and hence
$f(n,3,\ell)\leq (\ell+1)(n-\ell)+1$. \qed
\end{pf}

\begin{lem}{\upshape \cite{MaoLai}}\label{lem3-5}
Let $G$ be a connected graph of order $n\geq 6$. Let $k$ be an
integer with $3\leq k\leq n$. Then
$$
\tau_k(G)\leq \kappa(G)-k+2.
$$
\end{lem}

For $k=3$ and general $\ell \ (\frac{n-4}{3}\leq \ell\leq
\frac{n-r-2}{2})$, we give the upper and lower bound of
$f(n,3,\ell)$.
\begin{pro}\label{pro3-4}
Let $n,\ell$ be two integers with $\frac{n-4}{3}\leq \ell\leq
\frac{n-r-2}{2}$ and $n\geq 15$, where $n\equiv r \ (mod \ \ell+1)$
and $r\geq 2$. Then
$$
\left\lceil \frac{(\ell+2)n}{2}\right\rceil\leq f(n,3,\ell)\leq
\left\lfloor\frac{n}{\ell+1}\right\rfloor
(\ell+1)^2+(r+1)(\ell+1)+r-1.
$$
\end{pro}
\begin{pf}
Let $F=P_{s}\circ (\ell+1) K_{1}$ where
$s=\lfloor\frac{n}{\ell+1}\rfloor$. Set $P_{s}=u_1u_2\ldots u_{s}$.
Let $V((\ell+1) K_{1})=\{v_1,v_2\ldots v_{\ell+1}\}$. Let $G$ be the
graph obtained from $F$ by adding the vertices
$(u_{s+1},v_1),(u_{s+1},v_2),\ldots,(u_{s+1},v_r)$ and the edges
\begin{eqnarray*}
&& \{(u_{s},v_i)(u_{s+1},v_j)\,|\,1\leq i\leq \ell+1, \ 1\leq j\leq
r\}\cup \{(u_{s+1},v_j)(u_{s+1},v_{j+1})\,|\, 1\leq j\leq r-1\}\\
&&\cup \{(u_{1},v_j)(u_{1},v_{j+1})\,|\, 1\leq j\leq \ell\}.
\end{eqnarray*}
Clearly, $\kappa(G)=\ell+1$. From Lemma \ref{lem3-5}, $\tau_3(G)\leq
\kappa(G)-1=\ell$. We only need to show $\tau_3(G)\geq \ell$. Let
$H_1,H_2,\ldots, H_{s}$ be the copies corresponding to $(\ell+1)
K_{1}$, and let $H_{s+1}=\{(u_{s+1},v_j)\,|\,1\leq j\leq r\}$. For
any $S\subseteq V(G)$ and $|S|=3$, it suffices to show $\tau(S)\geq
\ell$. Let $S=\{x,y,z\}$.

Suppose that there exists some $H_{i} \ (1\leq i\leq s+1)$ such that
$|S\cap V(H_{i})|=3$. Without loss of generality, let $|S\cap
V(H_{1})|=3$. Then the trees induced by the edges in
$\{x(u_2,v_i),y(u_2,v_i),z(u_2,v_i)\}$ $(1\leq i\leq \ell+1)$ are
$\ell+1$ internally disjoint pendant $S$-Steiner trees in $G$, and
hence $\tau(S)\geq \ell+1$.

Suppose that there exist $H_{i},H_{j}$ such that $|S\cap
V(H_{i})|=2$ and $|S\cap V(H_{j})|=1$. Suppose $|j-i|\geq 2$. Note
that the subgraph induced by the vertices in
$\{(u_r,v_s)\,|\,i+1\leq r\leq j, \ 1\leq s\leq \ell+1\}$ is
$(\ell+1)$-connected. From Lemma \ref{lem2-1}, there exits a
$(z,U)$-fan, where $U=\{(u_{i+1},v_s)\,|\,1\leq s\leq \ell+1\}$, and
hence there is a path connecting $z$ and $(u_{i+1},v_s)$, say $P_s$,
where $1\leq s\leq \ell+1$. Then the trees induced by the edges in
$E(P_s)\cup \{(u_{i+1},v_s)x,(u_{i+1},v_s)y\}$ are $\ell+1$
internally disjoint pendant $S$-Steiner trees, and hence
$\tau(S)\geq \ell+1$. Suppose $|j-i|=1$. Without loss of generality,
let $S\cap V(H_{1})=\{x,y\}$ and $|S\cap V(H_{2})|=\{z\}$. Let
$x',y'$ be the vertices corresponding to $x,y$ in $H_2$, $z'$ be the
vertex corresponding to $z$ in $H_1$. Suppose $z'\not\in \{x,y\}$.
Without loss of generality, let $\{x,y,z'\}=\{(u_1,v_i)\,\,1\leq
i\leq 3\}$. Then the tree induced by the edges in
$\{xx',yx',x'(u_3,v_1),z(u_3,v_1)\}$, the tree induced by the edges
in $\{xy',yy',y'z',z'z\}$ and the trees induced by the edges in
$\{x(u_2,v_j),y(u_2,v_j),(u_1,v_j)(u_2,v_j),(u_1,v_j)z\} \ (4\leq
j\leq \ell+1)$ are $\ell$ internally disjoint pendant $S$-Steiner
trees in $G$, and hence $\tau(S)\geq \ell$. Suppose $z'\in \{x,y\}$.
Without loss of generality, let $z'=y$,
$\{x,y\}=\{(u_1,v_i)\,\,1\leq i\leq 2\}$. Then the tree induced by
the edges in $\{xx',yx',x'(u_3,v_1),z(u_3,v_1)\}$ and the trees
induced by the edges in
$\{x(u_2,v_j),y(u_2,v_j),(u_1,v_j)(u_2,v_j),(u_1,v_j)z\}  \ (3\leq
j\leq \ell+1)$ are $\ell$ internally disjoint pendant $S$-Steiner
trees in $G$, and hence $\tau(S)\geq \ell$.

Suppose that there exist $H_{i},H_{j},H_k$ such that $|S\cap
V(H_{i})|=|S\cap V(H_{j})|=|S\cap V(H_{k})|=1$. Without loss of
generality, let $|S\cap V(H_{1})|=|S\cap V(H_{2})|=|S\cap
V(H_{3})|=1$, and $x=(u_1,v_1)$, $y=(u_2,v_1)$, $z=(u_3,v_1)$. Then
the trees induced by the edges in
$\{x(u_2,v_j),y(u_2,v_j),(u_1,v_j)(u_2,v_j),(u_1,v_j)z\}  \ (2\leq
j\leq \ell+1)$ are $\ell$ internally disjoint pendant $S$-Steiner
trees in $G$, and hence $\tau(S)\geq \ell$.

From the above argument, we conclude that $\tau(S)\geq \ell$ for any
$S\subseteq V(G)$ and $|S|=3$. From the arbitrariness of $S$, we
have $\tau_3(G)\geq \ell$, and hence $\tau_3(G)=\ell$. So
$f(n,3,\ell)\leq \lfloor\frac{n}{\ell+1}\rfloor
(\ell+1)^2+(r+1)(\ell+1)+r-1$. From Proposition \ref{pro2-1},
$f(n,3,\ell)\geq \lceil \frac{(\ell+2)n}{2}\rceil$, as desired.\qed
\end{pf}

The following result is from \cite{MaoL}.
\begin{lem}{\upshape \cite{MaoL}}\label{lem3-6}
Let $G$ be a connected graph of order $n$. Then $\tau_3(G)=n-5$ if
and only if $\bar{G}$ is a subgraph of one of the following graphs.
\begin{itemize}
\item $C_i\cup C_j\cup (n-i-j)K_1 \ \ (i=3,4, \ j=3,4)$;

\item $C_i\cup \lfloor \frac{n-i}{2} \rfloor K_2  \ \ (i=3,4)$;

\item $P_5\cup \lfloor \frac{n-5}{2} \rfloor K_2$;

\item $C_i\cup (n-i) K_1  \ \ (i=5,6,7)$.
\end{itemize}
\end{lem}

\noindent \textbf{Proof of Theorem \ref{th1-3}:} From $(3)$ of
Theorem \ref{th1-1}, we have $f(n,3,0)=n-1$. From $(4)$ of Theorem
\ref{th1-1}, we have $f(n,3,1)=\lceil\frac{3n}{2}\rceil$. From
Proposition \ref{pro3-2}, $f(n,3,2)=2n$ for $n=pq$ and $p,q\geq 3$;
$2n\leq f(n,3,2)\leq \frac{5n}{2}$ for $n=2p$; $2n\leq f(n,3,2)\leq
4n-16$ if $n$ is a prime number. From Proposition \ref{pro3-3},
$\lceil \frac{(\ell+2)n}{2}\rceil\leq f(n,3,\ell)\leq
(\ell+1)(n-\ell)+1$ for $3\leq \ell\leq \frac{n-4}{3}$; From
Proposition \ref{pro3-4}, $\lceil \frac{(\ell+2)n}{2}\rceil\leq
f(n,3,\ell)\leq \lfloor\frac{n}{\ell+1}\rfloor
(\ell+1)^2+(r+1)(\ell+1)+r-1$ for $\frac{n-1}{3}\leq \ell\leq
\frac{n-r-2}{2}$ where $n\equiv r \ (mod \ \ell)$ and $r\geq 2$.
From $(5)$ of Theorem \ref{th1-1}, we have $\lceil
\frac{(\ell+2)n}{2}\rceil\leq f(n,3,\ell)\leq (\ell+2)(n-\ell-2)$
for $\frac{n-1}{3}\leq \ell\leq \frac{n-r-2}{2}$ where $n\equiv r \
(mod \ \ell)$ and $r=0,1$, or $\frac{n-r-2}{2}\leq \ell\leq
\frac{n-4}{2}$. From $(6)$ of Theorem \ref{th1-1}, we have
$\left\lceil \frac{(\ell+2)n}{2}\right\rceil\leq f(n,3,\ell)\leq
(\ell+2)(n-\ell-2)+{\ell+2\choose 2}$ for $\frac{n-2}{2}\leq
\ell\leq n-6$. From Lemma \label{lem3-6}, we have
$f(n,3,n-5)={n\choose 2}-\lfloor\frac{n-4}{2}\rfloor$. From $(1)$
and $(2)$ of Theorem \ref{th1-1}, $f(n,3,n-4)={n\choose 2}-2$ and
$f(n,3,n-3)={n\choose 2}$.\qed

\end{document}